\documentclass[11pt]{amsart}
\usepackage[T1]{fontenc}
\usepackage{amsmath,amssymb,mathtools}
\usepackage[hidelinks]{hyperref}
\usepackage[margin=1in]{geometry}
\numberwithin{equation}{section}
\newtheorem{theorem}{Theorem}[section]
\newtheorem{proposition}[theorem]{Proposition}
\newtheorem{corollary}[theorem]{Corollary}
\theoremstyle{remark}
\newtheorem{remark}[theorem]{Remark}
\newcommand{\R}{\mathbb R}
\newcommand{\T}{\mathbb T}
\newcommand{\divg}{\operatorname{div}}
\newcommand{\Ric}{\operatorname{Ric}}
\newcommand{\Def}{\operatorname{Def}}
\newcommand{\curl}{\operatorname{curl}}
\newcommand{\dx}{\,dx}
\newcommand{\dt}{\,dt}
\newcommand{\acc}{\dot\xi}
\newcommand{\ip}[2]{\langle #1,#2\rangle}
\title[Geometry of Navier--Stokes blowup]{The Landau--Raychaudhuri Equation and the Geometry of Finite-Time Singularities in Navier--Stokes Flows}
\author{Sergey E. Stepanov}
\address{Department of Mathematics, Financial University under the Government of the Russian Federation, 49--55 Leningradsky Prospect, Moscow 125468, Russian Federation}
\address{All-Russian Institute for Scientific and Technical Information of the Russian Academy of Sciences (VINITI RAS), 20 Usievicha Street, Moscow 125190, Russian Federation}
\email{s.e.stepanov@mail.ru}
\date{November 25, 2026}
\begin{document}
\begin{abstract}
We investigate the geometry of the unit velocity direction in the finite-time blowup solution of the forced three-dimensional incompressible Navier--Stokes equations constructed by OpenAI in September 2026. Using the spatial Landau--Raychaudhuri equation and the preserved axial data of that construction, we obtain exact rates of divergence of transverse rotation and expansion at the singular point. The deformation-to-rotation ratio tends to zero there, and the divergence of the streamline curvature vector is asymptotically the negative of the squared transverse rotation. We express the pressure Laplacian in terms of these geometric invariants and establish its leading axial asymptotics. In an inner profile region we prove divergence of the transverse shear, while in the exact azimuthal exterior we establish a sharp local space-time integrability threshold for streamline curvature. Weighted geometric identities, valid without nonvanishing-speed assumptions, explain why these singularities are compatible with finite integrated viscous dissipation. The results are geometric consequences of the cited construction, not an independent proof of blowup or a new general regularity criterion.
\end{abstract}
\maketitle

\section{Introduction}
The finite-time blowup of an incompressible viscous flow is a question about both the magnitude and the geometry of its velocity. Geometric approaches to regularity have included the coherence of the vorticity direction \cite{CF} and conditions on the velocity direction \cite{Vasseur,Guo}. The speed--direction formulation of the equations was studied in \cite{Olshanskii}. In \cite{Stepanov}, the spatial Landau--Raychaudhuri identity was combined with the deformation-based energy balance to describe the shear, rotation, expansion and curvature of the congruence defined by the unit velocity field. This approach is based on the Riemannian identity investigated in \cite{SM}.

On September 8, 2026, OpenAI released a construction of finite-time blowup for the forced three-dimensional Navier--Stokes equations \cite{OpenAI}. Its Theorem 1.1 gives, for every positive viscosity, a force $f\in C_c^\infty(\R^3\times(0,\infty);\R^3)$ and a flow starting from rest with bounded kinetic energy and unbounded velocity as $t\uparrow1$. In particular, the force is smooth on $\R^3\times[0,1]$ (by extension by zero near $t=0$) and across the singular time, whereas the velocity and pressure are smooth for $0\le t<1$. Corollary 10.6 of that work gives the periodic counterpart. The construction involves an anisotropically contracting vortex, a corrected axisymmetric background, annular oscillations and a final spatial localization. In the present paper we take that construction and its stated profile properties as input; we do not reprove its existence theorem. The purpose of the present paper is to give a geometrization of this AI-generated blowup result by means of the spatial Landau--Raychaudhuri equation. More precisely, we translate the analytic properties of the constructed solution into invariant statements about the congruence of instantaneous streamlines defined by $\xi=u/|u|$. The Landau--Raychaudhuri equation provides the central geometric relation between the transverse shear, rotation and expansion of this congruence and the divergence of its curvature vector. Combined with the Navier--Stokes energy and pressure equations, it reveals geometric balances that are not explicit in the original blowup theorem. This geometrization does not replace the analytic construction; rather, it gives an intrinsic description of the particular singular flow and relates its pointwise singularities to finite integrated viscous dissipation.

Two features of this example are particularly useful. First, the higher-order background corrections preserve selected axis traces, and the annular corrections do not affect the axis. This makes it possible to calculate geometric invariants exactly there. Second, the solution possesses an exact azimuthal exterior on a local region, which permits an independent analysis of streamline curvature without using the oscillatory part of the construction. We keep the distinction between the axis, a fixed inner similarity region and the exterior throughout. The main conclusions are the asymptotic dominance of transverse rotation on the axis, the limiting Landau--Raychaudhuri balance $\divg\acc/|\omega|^2\to-1$, and the leading pressure asymptotics there. We also identify off-axis regions of divergent transverse shear and streamline curvature and derive their sharp local space-time integrability thresholds. The results concern the forced construction and are not assertions about all Navier--Stokes singularities.

\section{Geometric identities and energy bounds}
Let $u$ be a smooth divergence-free vector field on a Riemannian three-manifold $(M,g)$. Write $s=|u|$ and, on $M_+=\{s>0\}$, set $\xi=u/s$, $P=\mathrm{Id}-\xi\otimes\xi^\flat$, $H=\xi^\perp$ and $g_H=g|_H$. The transverse derivative $B(X,Y)=g(\nabla_X\xi,Y)$, for $X,Y\in H$, has the orthogonal decomposition
\begin{equation}\label{decomp}
 B=\sigma+\omega+\tfrac12\theta g_H,\qquad
 \theta=\divg\xi,\qquad \acc=\nabla_\xi\xi.
\end{equation}
Here $\sigma$ is symmetric and trace-free, $\omega$ is skew-symmetric, and the dot denotes covariant differentiation along the spatial unit direction, not time differentiation. All norms are induced by $g$; the norms of $B,\sigma,\omega$ are Hilbert--Schmidt norms on $H$. Put $\nabla^\perp s=P\nabla s$ and $D_u=\Def(u)=\tfrac12\mathcal L_u g$. The following identities hold in dimension three:
\begin{align}
 \xi(s)&=-s\theta,\label{inc}\\
 |D_u|^2&=s^2(|\sigma|^2+\tfrac32\theta^2)
   +\tfrac12|\nabla^\perp s+s\acc|^2,\label{def}\\
 \xi(\theta)+\tfrac12\theta^2+|\sigma|^2-|\omega|^2
       +\Ric(\xi,\xi)&=\divg\acc.\label{LR}
\end{align}
Our curvature convention is the one for which the round unit sphere has positive Ricci curvature. Indeed, $\divg(s\xi)=0$ gives \eqref{inc}. From $\nabla_Xu=X(s)\xi+s\nabla_X\xi$, the components of $D_u$ relative to $TM=\langle\xi\rangle\oplus H$ are $D_u(\xi,\xi)=-s\theta$, $D_u(\xi,X)=\tfrac12 g(\nabla^\perp s+s\acc,X)$ for $X\in H$, and $D_u|_{H\times H}=s(\sigma+\tfrac12\theta g_H)$. Taking the squared norm gives \eqref{def}. Finally, the Ricci identity $\divg(\nabla_\xi\xi)=\xi(\divg\xi)+\operatorname{tr}((\nabla\xi)^2)+\Ric(\xi,\xi)$, together with $\operatorname{tr}((\nabla\xi)^2)=\operatorname{tr}(B^2)=|\sigma|^2-|\omega|^2+\tfrac12\theta^2$, gives \eqref{LR}. In particular, the divergence version of \eqref{LR} reads
\begin{equation}\label{LRdiv}
 \divg(\acc-\theta\xi)=|\sigma|^2-|\omega|^2-\tfrac12\theta^2+\Ric(\xi,\xi).
\end{equation}

\begin{proposition}[Weighted identities]\label{prop:weighted}
On $M_+$ one has
\begin{align}
 |\nabla u|^2&=|\nabla s|^2+s^2\bigl(|\sigma|^2+|\omega|^2+\tfrac12\theta^2+|\acc|^2\bigr),\label{grad}\\
 \divg(\nabla_u u)&=s^2\bigl(|\sigma|^2-|\omega|^2+\tfrac32\theta^2\bigr)
       +2s\ip{\nabla^\perp s}{\acc}+\Ric(u,u).\label{weighted}
\end{align}
The right-hand sides, with the directional terms extended by zero on $\{s=0\}$, are locally integrable and the identities hold almost everywhere on $M$. If $M$ is closed, or if $M=\R^3$ and $u$ is compactly supported, then
\begin{align}\label{integrated}
 \int_{M_+}s^2|\omega|^2\,dV_g
 &=\int_{M_+}s^2(|\sigma|^2+\tfrac32\theta^2)\,dV_g\\
 &\quad+2\int_{M_+}s\ip{\nabla^\perp s}{\acc}\,dV_g
       +\int_M\Ric(u,u)\,dV_g.\notag
\end{align}
No restriction on the dimension of the zero set of $u$ or on the nondegeneracy of its zeros is required.
\end{proposition}
\begin{proof}
Differentiation of $u=s\xi$ yields $\nabla_Xu=X(s)\xi+s\nabla_X\xi$. The summands are orthogonal, and $|\nabla\xi|^2=|B|^2+|\acc|^2$. This proves \eqref{grad}. By \eqref{inc}, $\nabla_u u=s^2(\acc-\theta\xi)$. Taking its divergence and applying \eqref{LRdiv}, with $\xi(s)=-s\theta$, gives \eqref{weighted}. All weighted quadratic terms are bounded locally by a constant times $|\nabla u|^2$ (apart from the smooth Ricci term), since $s\nabla\xi=P\nabla u$ and $|\nabla s|\le |\nabla u|$ almost everywhere. The locally Lipschitz function $s$ and the smooth components of $u$ have vanishing weak derivatives almost everywhere on their respective zero sets; hence $\nabla s=\nabla u=0$ almost everywhere on $\{s=0\}$. At those points $\divg(\nabla_u u)=0$ as well. This proves the almost-everywhere identities. Integrating the smooth divergence on the left of \eqref{weighted} gives \eqref{integrated}, without excising any zeros.
\end{proof}

Consider the forced system with constant density normalized to one and deformation viscosity $\nu>0$:
\begin{equation}\label{NS}
 \partial_tu+\nabla_u u=-\nabla p+2\nu(\divg D_u)^\sharp+f,
 \qquad \divg u=0.
\end{equation}
Here $p$ is pressure divided by the physical density. On a closed manifold, or for compactly supported smooth fields in Euclidean space, integration by parts gives
\begin{equation}\label{energy}
 \tfrac12\tfrac{d}{dt}\|u\|_2^2+2\nu\|D_u\|_2^2=\int_M g(f,u)\,dV_g.
\end{equation}

\begin{corollary}\label{cor:energy}
Suppose that $u$ is smooth for $0\le t<T<\infty$, $\sup_{t<T}\|u(t)\|_2<\infty$ and $f\in L^1(0,T;L^2)$. On a closed $M$ or for compactly supported Euclidean solutions as above,
\begin{equation}\label{bound}
 \int_0^T\int_{M_+(t)}s^2\bigl(|\sigma|^2+|\omega|^2+\tfrac12\theta^2+|\acc|^2\bigr)\,dV_g\dt<\infty.
\end{equation}
The space-time integral of $|\nabla s|^2$ is also finite.
\end{corollary}
\begin{proof}
Equation \eqref{energy} bounds $\int_0^T\|D_u\|_2^2\dt$. The Korn inequality on a closed Riemannian manifold gives $\|\nabla u\|_2^2\le C_g(\|D_u\|_2^2+\|u\|_2^2)$. In the flat case integration by parts and incompressibility give the exact identity $\|\nabla u\|_2^2=2\|D_u\|_2^2$. Apply \eqref{grad}.
\end{proof}

\begin{proposition}[Pressure and a continuation obstruction]\label{prop:pressure}
In flat space the pressure of \eqref{NS} satisfies on $M_+$
\begin{equation}\label{pressure}
 \begin{split}
 \Delta p={}&s^2\bigl(|\omega|^2-|\sigma|^2-\tfrac32\theta^2\bigr)\\
 &-2s\ip{\nabla^\perp s}{\acc}+\divg f.
 \end{split}
\end{equation}
The equality extends almost everywhere using the convention of Proposition \ref{prop:weighted}. Moreover, for a smooth solution on $\R^3$ or $\T^3$ of bounded kinetic energy which cannot be continued smoothly beyond a finite time $T$ with smooth compactly supported (or smooth periodic) forcing, one has
\begin{equation}\label{serrin}
 \int_0^T\biggl(\int_{M_+(t)}
       (|\nabla^\perp s|^2+s^2\theta^2)\,dx\biggr)^2\dt=\infty.
\end{equation}
\end{proposition}
\begin{proof}
Take the divergence of \eqref{NS} in the flat case, where the viscous term has zero divergence, and use \eqref{weighted}. By \eqref{inc}, the spatial integral in \eqref{serrin} is $\|\nabla s\|_2^2$. If its time integral after squaring were finite, the Sobolev inequality, together with the energy bound in the periodic case, would imply $u\in L^4(0,T;L^6)$. The classical Prodi--Serrin continuation criterion \cite{Prodi,Serrin} (with smooth forcing) then contradicts the assumed breakdown. This is a geometric rewriting of a classical criterion, not an independent regularity theorem.
\end{proof}

\section{The axial geometry of the September construction}
We use the viscosity-one solution in \cite{OpenAI}. The viscosity rescaling is discussed below. In order to avoid confusing the expansion $\theta$ with the cylindrical angle, denote the latter by $\varphi$. Let $\tau=1-t$ and use the exponents $a=\tfrac12+h$, $d=\tfrac12-h$, with $0<h<10^{-2}$. The similarity variables in \cite[Eq.~(3.2)]{OpenAI} are
\begin{equation}\label{coordinates}
 \tau=q(1-\eta^2),\qquad z=q^d\eta,\qquad X=\frac{r^2}{2q}.
\end{equation}
In particular, $q=\tau$ and $\eta=0$ at $z=0$. The leading azimuthal and axial velocities are $q^{-a}E(X,\eta)$ and $q^{-a}U(X,\eta)$, where $E=\sqrt{2X}\phi/C$ with $C>0$ fixed \cite[Eq.~(4.3)]{OpenAI}. The axis profiles satisfy
\begin{equation}\label{axisdata}
 U(0,\eta)=4\eta+j_0,\qquad \phi(0,\eta)=\phi_*(\eta)>0,
 \qquad \phi_*(0)=1,\quad j_0>0
\end{equation}
by \cite[Eqs.~(B.1), (B.3)]{OpenAI}.

It is essential that these are not merely leading-order axis data. In \cite[Lemma 5.1]{OpenAI} all positive-order profiles for $U$ and $\phi$ have zero axis traces. The cutoff summation of Proposition 5.5 preserves these traces, as does the annular correction procedure (see the proof of Proposition 9.9 and Section 10.1). The final localization is identical to the local solution near the origin for all times sufficiently close to one, by Proposition 10.1. Consequently, for these times the completed solution is axisymmetric near the axis and satisfies there
\begin{equation}\label{traces}
 w(z,t):=u_z(0,z,t)=q^{-a}(4\eta+j_0),\qquad
 b(z,t):=\lim_{r\downarrow0}\frac{u_\varphi}{r}
 =C^{-1}q^{-a-1/2}\phi_*(\eta).
\end{equation}
The neighborhood in which the annular corrections vanish may shrink as $t\uparrow1$; only the axis traces are used in the next theorem.

\begin{theorem}[Exact axial asymptotics]\label{thm:axis}
For the solution just specified, at $x=0$ and for $t$ sufficiently close to one,
\begin{align}
 s&=j_0\tau^{-1/2-h},&
 \theta&=-4j_0^{-1}\tau^{-1/2+h},\label{axis1}\\
 |\omega|^2&=2C^{-2}j_0^{-2}\tau^{-1},&
 \sigma&=0,\quad \acc=0.\label{axis2}
\end{align}
Furthermore, if $W_u=\operatorname{skew}\nabla u$ is the skew part of the physical velocity gradient, then
\begin{equation}\label{ratio}
 |D_u|^2=24\tau^{-2},\qquad |W_u|^2=2C^{-2}\tau^{-2-2h},\qquad
 \frac{|D_u|^2}{|W_u|^2}=12 C^2\tau^{2h}.
\end{equation}
Thus the last ratio and $\theta^2/|\omega|^2=8C^2\tau^{2h}$ tend to zero.
\end{theorem}
\begin{proof}
Near the axis, a smooth axisymmetric field has the Cartesian form $u_1=\alpha x_1-\beta x_2$, $u_2=\beta x_1+\alpha x_2$, $u_3=\gamma$, with $\alpha,\beta,\gamma$ smooth functions of $r^2,z,t$. At $r=0$, $\gamma=w$, $\beta=b$ and $2\alpha=-w_z$ by incompressibility. Equation \eqref{coordinates} gives $q_z=0$, $\eta_z=\tau^{-d}$ at $z=0$, so \eqref{traces} yields $w=j_0\tau^{-a}$, $w_z=4\tau^{-1}$, and $b=C^{-1}\tau^{-1-h}$. With the convention that the columns of the gradient matrix correspond to derivative directions, this gives
\begin{equation}\label{matrix}
 \nabla u(t,0)=\begin{pmatrix}
 -2\tau^{-1}&-C^{-1}\tau^{-1-h}&0\\
 C^{-1}\tau^{-1-h}&-2\tau^{-1}&0\\
 0&0&4\tau^{-1}
 \end{pmatrix}.
\end{equation}
Since $w>0$, $\xi=e_z$ on the axis locally. The transverse matrix of $\nabla\xi$ is the transverse block in \eqref{matrix} divided by $w$, so its symmetric part is pure trace and its skew part has two entries of modulus $b/w$. Also $\nabla_\xi\xi=0$ at the point. These observations give \eqref{axis1}--\eqref{axis2}. Taking symmetric and skew parts of \eqref{matrix} proves \eqref{ratio}.
\end{proof}

\begin{theorem}[Singular Landau--Raychaudhuri balance and pressure]\label{thm:LR}
In the setting of Theorem \ref{thm:axis}, at $x=0$ one has
\begin{equation}\label{LRaxis}
 \divg\acc=-\frac{2}{C^2j_0^2}\tau^{-1}
 +\left(1+2h+\frac{24}{j_0^2}\right)\tau^{-1+2h}.
\end{equation}
In particular,
\begin{equation}\label{limit}
 \lim_{t\uparrow1}\frac{\divg\acc}{|\omega|^2}=-1.
\end{equation}
The pressure satisfies
\begin{equation}\label{paxis}
 \Delta p=\frac{2}{C^2}\tau^{-2-2h}-24\tau^{-2}+O(1),\qquad
 \lim_{t\uparrow1}\frac{\Delta p}{s^2|\omega|^2}=1.
\end{equation}
\end{theorem}
\begin{proof}
On an axis segment where $w>0$, $\theta=-w_z/w$ and $\xi(\theta)=\partial_z\theta$. At fixed $t<1$, the relation $q-z^2q^{2h}=\tau$ gives the Taylor coefficients $q_z(0)=0$, $q_{zz}(0)=2\tau^{2h}$, $\eta_z(0)=\tau^{-d}$ and $\eta_{zz}(0)=0$. Substitution in \eqref{traces} gives $w_{zz}/w=-2a\tau^{-2d}$ at $z=0$, hence
\[
 \xi(\theta)=\left(2a+\frac{16}{j_0^2}\right)\tau^{-2d}.
\]
Equation \eqref{LR}, with $\Ric=0$, $\sigma=0$, and \eqref{axis1}--\eqref{axis2} proves \eqref{LRaxis} and \eqref{limit}. In \eqref{pressure} the cross term vanishes on the axis since $\acc=0$. Inserting the axial values gives \eqref{paxis}, using boundedness of $\divg f$ at time one (Theorem 1.1 of \cite{OpenAI}).
\end{proof}
\begin{remark}
The normalized rotation is not the vorticity itself. In dimension three, on $M_+$, $|\omega|^2=\tfrac12\ip{\xi}{\curl\xi}^2=\ip{u}{\curl u}^2/(2s^4)$. In the present example the axial vorticity equals $2b e_z$, parallel to the velocity. The limiting relation \eqref{limit} concerns a divergence, not the curvature vector itself, which vanishes at the origin. The pointwise dominance of rotation over strain does not imply a corresponding global integral dominance or a causal mechanism for blowup.
\end{remark}

\section{Transverse shear in the inner region}
The vanishing of the shear on the axis does not imply its boundedness near the singular point. We next use the nontrivial radial dependence of the leading profiles to prove the opposite.

\begin{theorem}[Off-axis shear blowup]\label{thm:shear}
For the solution of \cite{OpenAI}, there is a compact rectangle $\mathcal Q\subset\{0<X<X_a,\ -1<\eta<1\}$ of positive area, in the inner region free of annular corrections, such that on the corresponding time-dependent spatial set $Q_\tau$ one has, uniformly as $\tau\downarrow0$,
\begin{equation}\label{shearbound}
 s\asymp\tau^{-1/2-h},\qquad |\sigma|\asymp\tau^{-1/2}.
\end{equation}
Here $X_a$ is the inner edge of the active annulus in the construction. In particular, the spatial supremum of $|\sigma|$ tends to infinity as $t\uparrow1$. For any $p>0$ and sufficiently small $\epsilon>0$,
\begin{equation}\label{shearp}
 \int_{1-\epsilon}^{1}\int_{Q_{1-t}}|\sigma|^p\dx\dt<\infty
 \quad\Longleftrightarrow\quad p<5-2h.
\end{equation}
\end{theorem}
\begin{proof}
Use $R=\sqrt{2X}$ and let $\eta_*=-j_0/4$. At this axial parameter, $U(0,\eta_*)=0$ by \eqref{axisdata}. We first show that $U_X(0,\eta_*)\ne0$. Set $L=1-2h\eta^2$, $d_0=1-\eta^2$, and $\Pi_0(\eta)=\Pi(0,\eta)$, where $\Pi$ is the leading pressure profile of \cite{OpenAI}. The stress-free inner equation \cite[Eq.~(4.13)]{OpenAI}, evaluated at $X=0$ and $\eta=\eta_*$, together with the axial momentum source in Eq.~(4.9), yields
\begin{equation}\label{ux}
 -2L(\eta_*)U_X(0,\eta_*)
 =-4d\eta_*-d_0(\eta_*)\Pi_0'(\eta_*)+4a\eta_*\Pi_0(\eta_*).
\end{equation}
In this calculation the terms proportional to $XU_X$, $X\Pi_X$ and $U$ vanish, and $U_\eta(0,\eta)=4$. The pressure datum satisfies $\Pi_0<0$ and $\eta\Pi_0'(\eta)>0$ for $\eta\ne0$ by Section B.1 of \cite{OpenAI}. Since $\eta_*<0$, all three terms on the right of \eqref{ux} are positive. Hence $U_X(0,\eta_*)<0$.

At fixed $\eta=\eta_*$, let $g_*=\phi_*(\eta_*)/C>0$ and $k=U_X(0,\eta_*)/2\ne0$. Analyticity at the axis implies
\[
 E=g_*R+O(R^3),\qquad U=kR^2+O(R^4).
\]
For $R>0$ sufficiently small define $b_0=E/(E^2+U^2)^{1/2}$ and $c_0=U/(E^2+U^2)^{1/2}$. The leading normalized field, in the common tangential velocity normalization, is $\xi_0=b_0 e_\varphi+c_0 e_z$. Put $v_0=c_0 e_\varphi-b_0 e_z$. In the transverse plane of $\xi_0$ use the orthonormal basis $(e_r,v_0)$.

We now justify taking this limit for the completed velocity, including its derivatives. On fixed compact inner profile sets with $R>0$, the formal expansion (5.1) and its cutoff realization (Proposition 5.5 and its proof, in particular (5.46)) give the leading azimuthal and axial profiles with $O(q^{2h})$ errors in fixed profile derivatives, while the common-normalization radial velocity is $O(q^h)$, with the same derivative control. The annular corrections vanish on the inner region (Proposition 9.9 and its proof, Step 5); the localization in Proposition 10.1 is identical to the local solution on these shrinking sets. The field $q^a u$ therefore converges to $E e_\varphi+U e_z$, uniformly with radial and angular profile derivatives. Its physical $z$ derivatives, after multiplication by $\sqrt q$, are $O(q^h)$, by the chain rule in (4.2) of \cite{OpenAI}. In particular, on sets where $E^2+U^2$ is bounded below, $\xi\to\xi_0$ and the rescaled physical gradient $\sqrt q\nabla\xi$ has a limit given by the radial derivative in $R$ and the cylindrical connection; the axial derivative contributes zero in this limit.

Take an orthonormal basis of $\xi^\perp$ converging to $(e_r,v_0)$, for instance by projecting $e_r$ and orthogonalizing. In this basis the limit of the rescaled off-diagonal shear component is
\begin{equation}\label{offdiag}
 \lim_{q\downarrow0}\sqrt q\,\sigma_{12}
 =\tfrac12\left(c_0\partial_R b_0-b_0\partial_R c_0-\frac{b_0c_0}{R}\right).
\end{equation}
Indeed, the first two terms are $\ip{\partial_R\xi_0}{v_0}$, while the last is the radial component of $\nabla_{v_0}\xi_0$ after rescaling, since $\nabla_{e_\varphi}e_\varphi=-e_r/r$. The off-diagonal component is unaffected by subtracting the trace. At $\eta_*$, $b_0=1+O(R^2)$ and $c_0=(k/g_*)R+O(R^3)$, so the right side tends to $-k/g_*\ne0$ as $R\downarrow0$. Choose a fixed sufficiently small positive $R_*$ in the inner region. By continuity, on a closed rectangle of positive area about $(R_*^2/2,\eta_*)$ the limit in \eqref{offdiag} is bounded away from zero. The full rescaled gradient is bounded there, proving \eqref{shearbound} since $q\asymp\tau$.

Finally, at fixed time the similarity map has $r\,dr=q\,dX$ at fixed $\eta$, and $dz/d\eta=q^d L/(1-\eta^2)$ by \eqref{coordinates}. Thus
\begin{equation}\label{jac}
 dx=q^{1+d}\frac{L}{1-\eta^2}\,dX\,d\eta\,d\varphi.
\end{equation}
On the chosen rectangle this is comparable to $\tau^{3/2-h}dX\,d\eta\,d\varphi$. The integral in \eqref{shearp} is consequently comparable to $\int_0^\epsilon\tau^{3/2-h-p/2}\,d\tau$, which converges precisely for $p<5-2h$.
\end{proof}

\section{Curvature and dissipation in the azimuthal exterior}
We now use the exact exterior of the completed solution, not a leading-order approximation. In \cite[Theorem 3.1(iii), Eq.~(3.5)]{OpenAI}, on $X\ge X_{\rm ext}$ in the local domain the velocity is $u=K(r,\tau)e_\varphi$, with
\begin{equation}\label{K}
 K(r,\tau)=r^{-1-2h}\mathcal H(\tau/r^2),
\end{equation}
where $\mathcal H$ is positive and smooth on $[0,\infty)$, and $K_r<0$ (Lemma A.6 and the proof of Proposition 9.9, Step 4). The same formula applies within a fixed spatial neighborhood of the origin for the final solution at times sufficiently close to one by Proposition 10.1 of \cite{OpenAI}.

\begin{theorem}[Curvature threshold and local dissipation]\label{thm:ext}
Choose $X_2>X_1>X_{\rm ext}$ and $0<\eta_c<1$. Let $\Omega_\tau$ denote the spatial region corresponding to $X_1\le X\le X_2$, $|\eta|\le\eta_c$, with the full azimuthal range. For all sufficiently small $\tau>0$, the completed solution satisfies on this region
\begin{equation}\label{curv}
 \xi=e_\varphi,\qquad |\acc|^2=r^{-2},\qquad \sigma=\omega=\theta=0.
\end{equation}
Moreover,
\begin{equation}\label{scales}
 \begin{split}
 \operatorname{Vol}(\Omega_\tau)&\asymp\tau^{3/2-h},\qquad
 s^2\asymp\tau^{-1-2h},\qquad |D_u|^2\asymp\tau^{-2-2h},\\
 \int_{\Omega_\tau}s^2\dx&\asymp\tau^{1/2-3h},\qquad
 \int_{\Omega_\tau}|D_u|^2\dx\asymp\tau^{-1/2-3h}.
 \end{split}
\end{equation}
For $p>0$ and sufficiently small $\epsilon>0$ the curvature integral satisfies
\begin{equation}\label{critical}
 \int_{1-\epsilon}^1\int_{\Omega_{1-t}}|\acc|^p\dx\dt<\infty
 \quad\Longleftrightarrow\quad p<5-2h.
\end{equation}
In particular, the curvature does not belong to $L^5$ in space-time on the nonzero set of the velocity. The dissipation integral over these regions remains finite.
\end{theorem}
\begin{proof}
Positivity of $K$ gives $\xi=e_\varphi$. The cylindrical connection yields $\acc=-e_r/r$ and $\nabla_X\xi=0$ for $X\in\operatorname{span}(e_r,e_z)=\xi^\perp$, proving \eqref{curv}. On $\Omega_\tau$, $q\asymp\tau$ and $r\asymp\sqrt\tau$. Equation \eqref{jac} gives its volume. The ratio $\tau/r^2=(1-\eta^2)/(2X)$ ranges over a fixed compact interval, so \eqref{K} implies the asserted speed estimate. By \eqref{def} or direct differentiation,
\begin{equation}\label{Dext}
 |D_u|^2=\tfrac12(K_r-K/r)^2.
\end{equation}
As $K_r<0$ and $K>0$, $|K_r-K/r|\ge K/r$. Smoothness of the heat factor on the relevant compact interval and differentiation of \eqref{K} give the matching upper bound of order $r^{-2-2h}$. This proves \eqref{scales}. Its last exponent $-1/2-3h$ is greater than $-1$, so the time integral is finite. Finally, $|\acc|\asymp\tau^{-1/2}$ and \eqref{jac} reduce \eqref{critical} to the same power integral as in Theorem \ref{thm:shear}.
\end{proof}

\begin{corollary}\label{all}
In the September construction, each of $|\sigma|$, $|\omega|$, $|\theta|$ and $|\acc|$ is unbounded along a sequence of space-time points converging to $(0,1)$. In addition, both $\sigma$ and $\acc$ fail to belong to $L^5$ on the space-time nonzero set of the velocity. Nevertheless, all their speed-weighted squared integrals in \eqref{bound} are finite.
\end{corollary}
\begin{proof}
The axial statements follow from Theorem \ref{thm:axis}, the off-axis ones from Theorems \ref{thm:shear} and \ref{thm:ext}. The latter two give the $L^5$ obstructions since $h>0$. Bounded kinetic energy and compactly supported smooth forcing are supplied by Theorem 1.1 of \cite{OpenAI}; apply Corollary \ref{cor:energy} with $T=1$.
\end{proof}

\begin{remark}[Scope of the integrability statements]
The exponent $5-2h$ is a sharp threshold on the specified shrinking profile regions only. It is not asserted to be a sharp global threshold on the entire nonzero set, where additional singularities of the normalized direction could arise near velocity zeros or within the oscillatory annulus. Under the Navier--Stokes parabolic scaling $u_\lambda(x,t)=\lambda u(\lambda x,\lambda^2t)$, $\lambda>0$, the direction satisfies $\xi_\lambda(x,t)=\xi(\lambda x,\lambda^2t)$ and $\nabla\xi_\lambda(x,t)=\lambda(\nabla\xi)(\lambda x,\lambda^2t)$. For a space-time region $Q$ in the nonzero set, let $Q_\lambda=\{(x,t):(\lambda x,\lambda^2t)\in Q\}$. Then
\[
 \iint_{Q_\lambda}|\nabla\xi_\lambda|^p\dx\dt
 =\lambda^{p-5}\iint_Q|\nabla\xi|^p\dx\dt.
\]
Thus $L^5_{t,x}$ is the scale-invariant space for a first derivative of the direction in three spatial dimensions. Failure of the $L^5$ condition for this solution is not a proof that this condition alone suffices for regularity.
\end{remark}

\section{Viscosity rescaling and concluding remarks}
The viscosity-one construction in \cite{OpenAI} gives the solution at viscosity $\nu>0$ by
\[
 \begin{aligned}
 u_\nu(x,t)&=\sqrt\nu\,u(x/\sqrt\nu,t),& p_\nu(x,t)&=\nu p(x/\sqrt\nu,t),\\
 f_\nu(x,t)&=\sqrt\nu\,f(x/\sqrt\nu,t).&&
 \end{aligned}
\]
Hence $\xi_\nu(x,t)=\xi(x/\sqrt\nu,t)$. The first-order normalized invariants are multiplied by $\nu^{-1/2}$ and $\divg\acc$ by $\nu^{-1}$ at corresponding points; $\nabla u$ and $\Delta p$ are unchanged at corresponding points. All rates and dimensionless limits above are preserved. The spatial integrals acquire the appropriate constant Jacobian factor, leaving convergence thresholds unchanged. The periodic construction of \cite[Corollary 10.6]{OpenAI} likewise retains the local geometric statements under its spatial rescaling.

The results reveal a distinction between the kinematic singularities of the normalized velocity field and the global energy dissipation of the physical velocity. At the spatial origin the transverse rotation dominates the expansion and the physical skew gradient dominates the deformation tensor. The Landau--Raychaudhuri identity relates this rotation to the negative divergence of the curvature vector even though the curvature vector itself is zero there. Off the axis, nontrivial radial changes of the velocity direction produce shear blowup; in the exact exterior, the direction has circular integral curves and unbounded curvature but no transverse shear, rotation or expansion. The weighted identities in Section 2 describe why the integral energy law is compatible with these different pointwise singularities. Although \eqref{NS} is forced ($f\not\equiv0$), the spatial Landau--Raychaudhuri identity is purely kinematic and contains no explicit forcing term, while $\divg f$ remains bounded at the singular time. Thus the leading singular terms in the axial geometric and pressure balances come from the contracting, anisotropically self-similar velocity profile rather than from a singular external force. This does not imply that the same evolution occurs without forcing. These conclusions use only the particular profile data and correction properties cited above; they make no claim of universality for other forced or unforced flows.

A further problem is to develop continuation criteria that control the geometry of the velocity direction jointly with its magnitude, in spaces compatible with the Navier--Stokes scaling. The present example supplies explicit tests that such proposed criteria must pass, and indicates the need to distinguish directional curvature from rotation and to account for the possible zero set of the velocity.

\medskip
\noindent\textbf{Data availability.} No new datasets were generated. The cited construction is publicly available in \cite{OpenAI}.

\noindent\textbf{AI assistance.} ChatGPT (OpenAI) was used to assist with mathematical derivations, source checking and preparation of this manuscript.

\end{document}